\documentclass[11pt]
{amsart}
\usepackage{amsthm,amsmath,amssymb,amsfonts}
\usepackage[backref,colorlinks]{hyperref}
\usepackage[margin=1.4in]{geometry}

\newtheorem{theorem}{Theorem}
\newtheorem{thm}[theorem] {Theorem}
\newtheorem{lemma}[theorem]{Lemma}

\newtheorem{claim}[theorem]{Claim}

\newtheorem{defi}[theorem]{Definition}

\let\eps=\varepsilon
\let\theta=\vartheta
\let\rho=\varrho
\let\sigma=\varsigma

\let\polishlcross=\l
\def\l{\ifmmode\ell\else\polishlcross\fi}

\def\vG{{\vec{G}}}

\def\vC{{\vec{C}}}

\def\vH{{\vec{H}}}

\def\E{{\bf E}}

\def\d{\delta}

\def\D{\Delta}

\def\l{\lambda}

\def\Pr{\mbox{{\bf Pr}}}

\newcommand{\proofstart}{{\bf Proof\hspace{2em}}}

\newcommand{\proofend}{\hspace*{\fill}\mbox{$\Box$}\vspace*{7mm}}

\newcommand{\Exp}{\mbox{\bf E}}

\newcommand{\rt}{\right}

\newcommand{\lt}{\left}

\begin{document}

\title[On the number of Hamilton cycles in sparse random graphs]{On the number of Hamilton cycles in sparse random graphs}

\author[R.~Glebov]{Roman Glebov}
\address{Institut f\"{u}r Mathematik, Freie Universit\"at Berlin, Arnimallee 3-5, D-14195 Berlin, Germany}
\email{glebov@mi.fu-berlin.de}
\thanks{The first author was supported by DFG within the research training group "Methods for Discrete Structures".}
\author[M.~Krivelevich]{Michael Krivelevich}
\address{School of Mathematical Sciences,
Sackler Faculty of Exact Sciences,
Tel Aviv University,
Tel Aviv 69978,
Israel}
\email{krivelev@post.tau.ac.il}
\thanks{The second author was supported in part by a USA-Israel BSF grant and by a grant from Israel Science
Foundation.}

\date{\today}
\begin{abstract}

We prove that the number of Hamilton cycles in the random graph $G(n,p)$ is $n!p^n(1+o(1))^n$ a.a.s.,
provided that $p\geq \frac{\ln n+\ln\ln n+\omega(1)}{n}$.
Furthermore, we prove the hitting-time version of this statement,
showing that in the random graph process, the edge that creates a graph of minimum degree $2$
creates $\lt(\frac{\ln n}{e}\rt)^n(1+o(1))^n$ Hamilton cycles a.a.s.
\end{abstract}

\maketitle

\setcounter{footnote}{0}
\renewcommand{\thefootnote}{\fnsymbol{footnote}}

\section{Introduction}

The goal of this paper is to estimate the number of Hamilton cycles in the random graph $G(n,p)$.
To be more formal, we show that the number of Hamilton cycle is a.a.s. concentrated around the expectation up to a factor $(1+o(1))^n$,
provided the minimum degree is at least $2$.

It is well known (see e.g.~\cite{bollobasbook}) that the minimum degree of $G(n,p)$ is a.a.s. at most one for $p\leq \frac{\ln n+\ln\ln n-\omega(1)}{n}$,
and therefore $G(n,p)$ contains no Hamilton cycle in this range of $p$ a.a.s.
Koml\'os and Szemer\'edi~\cite {KS} and Korshunov~\cite{Korshunov} were the first to show that this bound is tight,
i.e., $G(n,p)$ is a.a.s. Hamiltonian for every $p\geq \frac{\ln n+\ln\ln n+\omega(1)}{n}$.
Bollob\'as~\cite{bollobas} and
independently Ajtai, Koml\'os, and Szemer\'edi~\cite {AKS} proved the hitting time version of the above statement,
showing that in the random graph process, the very edge that increases the minimum degree to two also makes the graph Hamiltonian a.a.s.

There exists a rich literature about hamiltonicity of $G(n,p)$ in
the range when it is a.a.s. Hamiltonian. Recent results include
packing and covering problems (see e.g.~\cite{FK1}, \cite{KKO1},
\cite{KKO2}, \cite{MichaelWojciech}, \cite{us+T}, and~\cite{HKLO}),
local resilience (see e.g.~\cite{resilience}, \cite{FK2},
\cite{BSKS} and~\cite{LS}) and Maker-Breaker games (\cite{SS05},
\cite{HKSS09}, \cite{BSFHK12}, and~\cite{FGNK}). In this paper, we
are interested in estimating the typical number of Hamilton cycles
in a random graph when it is a.a.s. Hamiltonian. Several recent
results about Hamiltonicity (\cite{KKO2},\cite{MichaelWojciech},
\cite{resilience},\cite{FGNK}) can be used to show fairly easily that
$G(n,p)$ with $p=p(n)$ above the threshold for Hamiltonicity
contains typically many, or even exponentially many Hamilton cycles.
Here we aim however for (relatively) accurate bounds.

Using linearity of expectation we immediately see that the expected value of the number of Hamilton cycles in $G(n,p)$ is $\frac{(n-1)!}{2}p^n$.
As the common intuition for random graphs may suggest, we expect
the random variable to be concentrated around its mean, perhaps after some normalization (it is easy
to see that the above expressions for the expectation become exponentially large in $n$ already for $p$
inverse linear in $n$).

The reality appears to confirm this intuition -- to a certain
extent. Denoting by $X$ the number of Hamilton cycles in $G(n,p)$,
we immediately obtain $X<\left(\frac{np}{e}\right)^n$ a.a.s. by
Markov's inequality. Janson~\cite{Janson}
considered the distribution of $X$ for $p=\Omega\lt(\sqrt{n}\rt)$ and proved that $X$ is log-normal distributed,
implying that $X=\lt(\frac{np}{e}\rt)^n(1+o(1))^n$ a.a.s.
It is instructive to observe that assuming $p=o(1)$, the distribution
of $X$ is in fact concentrated way below its expectation, in particular implying that $X/\E(X)\overset{p}{\rightarrow}0$.
For random graphs of density $p =o\lt( n^{-1/2}\rt)$ not much appears to be known about the asymptotic behavior
of the number of Hamilton cycles in corresponding random graphs.
We nevertheless mention the result
of Cooper and Frieze~\cite{CF}, who proved that in the random graph process typically at the very moment the
minimum degree becomes two, not only the graph is Hamiltonian but it has $(\log n)^{(1-o(1))n}$ Hamilton
cycles.

Our main result is the following theorem, which can be interpreted
as an extension of Janson's results~\cite{Janson} to the full range
of $p(n)$.

\begin{thm}
\label{main}
Let $G\sim G(n,p)$ with $p\geq \frac{\ln n+\ln\ln n+\omega(1)}{n}$. Then the number of Hamilton cycles is
$n!p^n(1-o(1))^n$
a.a.s.
\end{thm}

Improving the main result of~\cite{CF}, we also show the following
statement.

\begin{thm}
\label{hittingtime}
In the random graph process, at the very moment the minimum degree becomes two,
the number of Hamilton cycles becomes $(\ln n/e)^n(1-o(1))^n$ a.a.s.
\end{thm}

We continue with a short overview of related results for other models of random and pseudorandom graphs.
For the model $G(n,M)$ of random graphs with $n$ vertices and $M$ edges, notice the result of Janson~\cite{Janson}
showing in particular that for the regime $n^{3/2}\ll M \leq 0.99\binom{n^2}{m}$,
the number of Hamilton cycles is indeed concentrated around its expectation.
The situation appears to change around
$M = \Theta\lt(n^{3/2}\rt)$, where the asymptotic distribution becomes log-normal instead.
Notice also that the number
of Hamilton cycles is more concentrated in $G(n,M)$ compared to $G(n,p)$; this is not surprising as
$G(n,M)$ is obtained from $G(n, p)$ by conditioning on the number of edges of $G$ being exactly equal to
$M$, resulting in reducing the variance.

For the probability space of random regular graphs, it is the opposite case of very sparse graphs that
is relatively well understood. Janson~\cite{Jansonreg}, following the previous work of Robinson and Wormald~\cite{RW1},~\cite{RW2},
described the asymptotic distribution of the number of Hamilton cycles in a random $d$-regular
graph $G(n,d)$ for a {\em constant} $d \geq 3$. The expression obtained is quite complicated, and we will not
reproduce it here.
For the case of growing degree $d = d(n)$,
the result of Krivelevich~\cite{Michael} on the number of Hamilton cycles in $(n,d,\lambda)$-graphs in addition to known eigenvalue
results for $G_{n,d}$ imply an estimation on the number of Hamilton cycles in $G_{n,d}$ with a superpolylogarithmic lower bound on $d$.

For an
overview of these results as well as of the corresponding results in pseudorandom settings,
we refer the interested reader to~\cite{Michael}.

\subsection{Definitions and notation}

The random oriented graph $\vG(n,p)$ is obtained from $G(n,p)$ by randomly giving an orientation to every edge
(every of the two possible directions with probability $1/2$).
Notice that whenever we use the notation $\vG$ for an oriented graph, there exists an underlying non-oriented graph
obtained by omitting the orientations of the edges of $\vG$; it is denoted by $G$.
Making the notation consistent, when omitting the vector arrow above an oriented graph,
we refer to the underlying non-oriented graph.

Given a graph $G$, we denote by $h(G)$ the number of Hamilton cycles in $G$.
We call a spanning $2$-regular subgraph of $G$ a $2$-{\em factor}.
Notice that every connected component of a $2$-factor is a cycle.
We denote by $f(G,s)$ the number of $2$-factors in $G$ with exactly $s$ cycles.
Similarly, a $1$-factor of an oriented graph $\vG$ is a spanning $1$-regular subgraph,
i.e., a spanning subgraph with all in- and outdegrees being exactly one.
Analogously, the number of $1$-factors in $\vG$ with exactly $s$ cycles is denoted by $f\lt(\vG,s\rt)$.
For the purposes of our proofs, we relax the notion of a $2$-factor and call a spanning subgraph $H\subseteq G$ an
{\em almost $2$-factor} of $G$ if $H$ is a collection of vertex-disjoint cycles
and at most $|V(G)|/\ln^2(|V(G)|)$ isolated vertices.
We denote the number of almost $2$-factors of $G$ containing exactly $s$ cycles by $f'(G,s)$.
Similarly to the notation for non-oriented graphs, we call an oriented subgraph $\vH$ of $\vG$ an
{\em almost $1$-factor} of $\vG$ if $\vH$ is a $1$-regular oriented graph
on at least $|V(\vG)|-|V(\vG)|/\ln^2(|V(\vG)|)$ vertices.
The number of almost $1$-factors of $\vG$ with exactly $s$ cycles is denoted by $f'\lt(\vG,s\rt)$.

As usual, in a graph $G$ for a vertex $x\in V(G)$ we denote by $d_G(x):=|N_G(x)|$ its {\em degree}, i.e., the size of its neighborhood.
We denote by $\d(G)$ and respectively $\D(G)$ its minimum and maximum degrees.
For a set $S\subseteq V(G)$, we denote by $N_G(S)$ the set of all vertices outside $S$ having a neighbor in $S$.
Whenever the underlying graph is clear from the context we might omit the graph from the index.
Similarly, in an oriented graph $\vG$ for a vertex $x\in V\lt(\vG\rt)$ we call $d_{in,\vG}(x):=\lt| \lt\{y\in V\lt(\vG\rt):\, yx \in E\lt(\vG\rt)\rt\}\rt|$
the {\em indegree} of $x$
and $d_{out,\vG}(x):=\lt| \lt\{z\in V\lt(\vG\rt):\, xz \in E\lt(\vG\rt)\rt\}\rt|$ the {\em outdegree} of $x$.
We denote by $\d_{in}\lt(\vG\rt)$, $\D_{in}\lt(\vG\rt)$, $\d_{out}\lt(\vG\rt)$, and $\D_{out}\lt(\vG\rt)$ the minimum and maximum in- and outdegrees of $\vG$.

In a graph $G$ for two sets $A,B\subseteq V(G)$ we denote by $e_G(A,B)$ the number of edges incident with both sets.
In an oriented graph $\vG$, for two sets $A,B\subseteq V(G)$
the notation $e_\vG(A,B)$ stands for the number of edges going from a vertex in $A$ to a vertex in $B$.
We write $e_G(A):=e_G(A,A)$ and $e_\vG(A):=e_\vG(A,A)$ for short.
Similarly to the degrees, whenever the underlying graph is clear from the context we might omit the graph from the index.

To simplify the presentation, we omit all floor and ceiling signs whenever these are not crucial.
Whenever we have a graph on $n$ vertices, we suppose its vertex set to be $[n]$.

\subsection{Outline of the proofs}

In Section~\ref{proofs}, the lower bounds for Theorems~\ref{main} and~\ref{hittingtime} are proven in the following steps.
\begin{itemize}
\item
In Lemma~\ref{orient} we show using the permanent of the incidence matrix that under certain pseudorandom conditions,
an oriented graph contains sufficiently many oriented $1$-factors.
\item
In Lemma~\ref{core} we prove that the random oriented graph $\vG(n,p)$ a.a.s. contains a large subgraph
with all in- and outdegrees being concentrated around the expected value.
This subgraph then satisfies one of the conditions of Lemma~\ref{orient}.
\item
In Lemma~\ref{number2-factors} we show that
the random graph $G(n,p)$ contains many almost $2$-factors a.a.s.
In the proof, we orient the edges of $G(n,p)$ randomly and apply Lemma~\ref{orient} to the subgraph with almost equal degrees
whose existence is guaranteed by Lemma~\ref{core} a.a.s.
\item
In Lemma~\ref{shortcycles} we prove that most of these almost $2$-factors have few cycles a.a.s.
\item
We then call a graph $p$-expander if it satisfies certain expansion properties
and show in Lemma~\ref{expander} that in the random graph process, the graph $G(n,p)$ has these properties in a strong way.
\item
Lemma~\ref{St} shows that in any graph having the $p$-expander properties and minimum degree $2$,
for any path $P_0$ and its endpoint $v_1$
many other endpoints can be created by a small number of rotations with fixed endpoint $v_1$.
\item
Lemma~\ref{digest} contains the main technical statement of the paper.
Its states that in a graph satisfying certain pseudorandom conditions,
for almost every almost $2$-factor $F$ with few components, there exists a Hamilton cycle with a small Hamming distance from $F$.
The proof is a straightforward use of Lemma~\ref{St}.
\item
The proofs of Theorems~\ref{main} and~\ref{hittingtime} are
completed with a double counting argument. On the one hand, by
Lemma~\ref{shortcycles} there exist many almost $2$-factors with few
cycles a.a.s. Furthermore, for each of these almost $2$-factors
there exists a Hamilton cycle with small Hamming distance from it
a.a.s. by Lemma~\ref{digest}. On the other hand, for each Hamilton
cycle, there are not many almost $2$-factors with few cycles having
a small Hamming distance from it. Hence, the number of Hamilton
cycles is strongly related to the number of almost $2$-factors with
few cycles, finishing the proof.
\end{itemize}

\section{The proofs}
\label{proofs}

Let $G\sim G(n,p)$.
Since $\E(h(G))=(n-1)!p^n/2$, we obtain \[h(G)<\ln n (n-1)!p^n/2<\left(\frac{np}{e}\right)^n\] a.a.s., using just Markov's inequality.
Thus, for the remainder of the section we are only interested in the lower bound on the typical number of Hamilton cycles in the random graph.

We know from~\cite{KKO1} and using e.g. the results from~\cite{KKO2} and~\cite{MichaelWojciech} that in $G\sim G(n,p)$ there are at least
$\lt(\frac{\lt\lfloor\delta(G)/2\rt\rfloor}{n}\rt)^n n!$ $2$-factors a.a.s.
We now want to give an a.a.s. lower bound on the number of $2$-factors in $G$, and we want to do it
within a multiplicative error term of at most $2^{o(n)}$ from the ``truth'', basically deleting the $2$ from the denominator in the above expression in the case $p\gg \ln n/n$,
and replacing the term $\lt\lfloor\delta(G)/2\rt\rfloor$ by asymptotically $np$.

We first prove a pseudo-random technical statement that will give us the desired inequality once we show that $G$
(or a large subgraph of it) satisfies the pseudo-random conditions.
The proof is based on the permanent method as used in~\cite{FK}.

\begin{lemma}
\label{orient}
Let $n$ be sufficiently large
and $r\gg \ln\ln n$,
and let $\vec{G}$ be an oriented graph on $n$ vertices satisfying the following (pseudo-random) conditions:
\begin{itemize}
\item
$\delta_{in}(\vec{G}), \delta_{out}(\vec{G}),
\Delta_{in}(\vec{G}), \Delta_{out}(\vec{G})
\in \lt(r-4r/\ln\ln n, r+4r/\ln\ln n\rt)$
\item
for any two sets $A,B\subset V(\vec{G})$ of size at most $|A|,|B|\leq 0.6n$,
there are at most $0.8r\sqrt{|A||B|}$ edges going from $A$ to $B$.
\end{itemize}
Then $\vec{G}$ contains at least $\lt(\frac{r-100r/\ln\ln n}{e}\rt)^n$ oriented $1$-factors.
\end{lemma}

\proofstart
Create an auxiliary bipartite graph $G'$ from $\vec{G}$ in the following way:
take two copies $X$ and $Y$ of the vertex set $[n]$
by doubling each vertex $v\in [n]$ into $v_X\in X$ and $v_Y\in Y$.
We put a (non-oriented) edge $uv\in E(G')$ between vertices $u_X\in X$ and $v_Y \in Y$
if $\vec{uv}\in E(\vec{G})$ is an edge oriented from $u$ to $v$ in $\vec{G}$.
We observe a one-to-one correspondence between oriented $1$-factors in $\vec{G}$
and perfect matchings in $G'$.

In order to use the permanent to obtain a lower bound on the number of perfect matchings of $G'$,
we need a (large) spanning regular subgraph of $G'$. Its existence is guaranteed by the following claim.

\begin{claim}
\label{oreryser}
$G'$ contains a spanning regular subgraph $G''$ with regularity at least $d=r-100r/\ln\ln n$.
\end{claim}

\proofstart
Applying the Ore-Ryser theorem~\cite{oreryser} we see that the statement of the claim is true provided that for every $Y'\subseteq Y$ we have
\[ d|Y'|\leq \sum_{x\in X}\min \{d, e_{G'}(x,Y')\}.\]
Suppose to the contrary that this conidition does not hold, i.e., there exists a $Y'\subseteq Y$ s.t.
\[ d|Y'|> \sum_{x\in X}\min \{d, e_{G'}(x,Y')\}.\]
We examine the number of edges incident to $Y'$ that can be deleted from $G'$ without disturbing the right hand side
of the above inequality. Formally, we denote it by $c= \sum_{x\in X}\max \{0, e_{G'}(x,Y')-d\}$.
Notice that
\begin{align*}
 c&= e(X, Y')-\sum_{x\in X}\min \{d, e_{G'}(x,Y')\} > e(X, Y')-d|Y'|
\end{align*}
as supposed above.

Since $(r-4r/\ln\ln n)|Y'|\leq \delta(G')|Y'|\leq e(X,Y')<d|Y'|+c$, we obtain
\[c>\frac{96r}{\ln\ln n}|Y'|.\]
On the other hand, denoting by $X'$ the set of vertices that have at least $d$ neighbors in $Y'$,
and noticing that $\Delta(G')\leq r+4r/\ln\ln n$,
we obtain
\[c\leq\frac{104r}{\ln\ln n}|X'|.\]
Hence,
\begin{equation}
\label{x'>}
|X'|> 0.9|Y'|.\end{equation}

Notice that by the choices of $Y'$ and $X'$, we have
\begin{equation}
d|Y'|> \sum_{x\in X}\min \{d, e_{G'}(x,Y')\}=d|X'|+e(X\setminus X',Y'). \label{mark1}
\end{equation}

For the number of edges between $Y\setminus Y'$ and $X\setminus X'$ we see that
\begin{align}
(r+4r/\ln\ln n)|Y\setminus Y'|&\geq e(X\setminus X',Y\setminus Y')= e(X\setminus X',Y)-e(X\setminus X',Y')\nonumber\\
&\overset{\eqref{mark1}}{>} \delta(G')|X\setminus X'|-d(|Y'|-|X'|)\nonumber\\
&\geq\frac{96r}{\ln\ln n} |X\setminus X'|
+(r-100r/\ln\ln n)|Y\setminus Y'|,
\end{align}
leading to
\begin{equation}
\label{x'<}
|X\setminus X'|< 1.1|Y\setminus Y'|.\end{equation}

Furthermore, notice that by~\eqref{mark1} it holds that
\begin{equation}
\label{x'=}
|X'|< |Y'|.\end{equation}

We prove the claim by case analysis.
\begin{itemize}
\item
If $|Y'|\leq n/2$, we obtain for the number of edges between $X'$ and $Y'$
\[e(X',Y')\overset{\mbox{Choice of }X'}{\geq} d|X'|\overset{\eqref{x'>}}{>} 0.9d\sqrt{|X'||Y'|}>0.8r\sqrt{|X'||Y'|},\]
contradicting the second condition of the lemma.
\item
If $|Y'|>n/2$, then again by the definition of $X'$ we obtain $e(X',Y')\geq d|X'|$,
leading to
\[e(X\setminus X', Y')\overset{\eqref{mark1}}{<} d(|Y'|-|X'|)=d(|X\setminus X'|-|Y\setminus Y'|)
\overset{\eqref{x'<}}{<}0.1d|Y\setminus Y'|\overset{\eqref{x'=}}{<}0.1d|X\setminus X'|.\]
Thus, using the fact that $\delta(G')>d$, we see that
\[e(X\setminus X',Y\setminus Y')\geq 0.9d|X\setminus X'|
\overset{\eqref{x'=}}{>}0.8r\sqrt{|X\setminus X'|\cdot |Y\setminus Y'|},\]
again contradicting the same condition of the lemma,
since now both $X\setminus X'$ and $Y\setminus Y'$ have size less than $0.6n$ by~\eqref{x'<}.
\end{itemize}

\proofend

We observe that the number of perfect matchings in $G''$ equals the permanent of the incidence matrix of $G''$.
Hence the result of
Egorychev~\cite{Egorychev} and Falikman~\cite{Falikman} on the conjecture of van der Waerden implies
that the number of perfect matchings in $G''$ is at least $d^nn!/n^n>\lt(\frac{d}{e}\rt)^n$.
\proofend

In order to use Lemma~\ref{orient},
we first prove the a.a.s. existence of a large subgraph of $\vG(n,p)$
satisfying the degree-conditions of Lemma~\ref{orient} a.a.s.

\begin{lemma}
\label{core}
Let $\vG\sim \vG(n,p)$ with $p\geq \ln n/n$.
Then there exists a set $V'\subseteq[n]$ of at least $n-n/\ln^2n$ vertices of $\vG$
such that the graph $\vC:=\vG[V']$ satisfies
$\delta_{in}(\vec{C}), \delta_{out}(\vec{C}),
\Delta_{in}(\vec{C}), \Delta_{out}(\vec{C})
\in \lt(\frac{np-3np/\ln\ln n}{2}, \frac{np+np/\ln\ln n}{2}\rt)$
a.a.s.
\end{lemma}

\proofstart
We observe using Chernoff's inequality (see e.g. Corollary~A.1.14 in~\cite{alonspencer})
that for $p\gg \ln n(\ln\ln n)^2/n$ the statement holds for $V'=[n]$ a.a.s.
Hence, from now on we assume $np=O(\ln n(\ln\ln n)^2)$.

Let $L$ be the set of all vertices whose in- or outdegree is at most $\frac{np-np/\ln\ln n}{2}+1$.
For every $y\in [n]$, we can estimate
using Chernoff's inequality (see e.g. Corollary~A.1.14 in~\cite{alonspencer})
\[\Pr(y\in L)=\exp\lt(-\Omega\lt(\ln n /(\ln\ln n)^2\rt)\rt). \]
Thus, by Markov's inequality we obtain
\begin{align}
|L|&\leq \ln n\cdot\E(|L|)=\ln n \cdot(n-1)\exp\lt(-\Omega\lt(\ln n /(\ln\ln n)^2\rt)\rt)\nonumber\\
&= n\exp\lt(-\Omega\lt(\ln n /(\ln\ln n)^2\rt)\rt)\label{|L_x|}
\end{align}
a.a.s.

Fix an arbitrary vertex $x\in [n]$.
We denote
\[L_x=\lt\{y\in [n]\setminus \{x\}:\,
d_{in,\vG-x}(y)\leq\frac{np-np/\ln\ln n}{2}
\mbox{ or }d_{out,\vG-x}(y)\leq\frac{np-np/\ln\ln n}{2}\rt\}.\]
Notice that $L_x\subseteq L$, and thus~\eqref{|L_x|} bounds $|L_x|$ as well.

Since for every $y\in [n]\setminus \{x\}$ the events ``$xy \in E(G)$''
and ``$y\in L_x$'' are independent,
we obtain using Chernoff's inequality again (see e.g. Theorem~A.1.12 in~\cite{alonspencer})
\begin{align}
\label{n(x)capl}
&\Pr \lt[\lt(|N_G(x)\cap L_x|\geq \frac{np}{2\ln\ln n}\rt)\,|\, \lt(|L_x|= n\exp\lt(-\Omega\lt(\ln n /(\ln\ln n)^2\rt)\rt)\rt)\rt]\nonumber\\
\qquad &\leq \exp\lt(-\frac{np}{2\ln\ln n}\Omega\lt(\ln n /(\ln\ln n)^2\rt)\rt)=o(1/n).
\end{align}

Similarly, we let
$R$ be the set of all vertices whose in- or outdegree is at least $\frac{np-np/\ln\ln n}{2}-1$
and obtain
\begin{equation}
\label{|R_x|}
|R|\leq  n\exp\lt(-\Omega\lt(\ln n /(\ln\ln n)^2\rt)\rt)
\end{equation}
a.a.s.

We define analogously
\[R_x=\lt\{y\in [n]\setminus \{x\}:\,
d_{in,\vG-x}(y)\geq\frac{np+np/\ln\ln n}{2}-1
\mbox{ or }d_{out,\vG-x}(y)\geq\frac{np+np/\ln\ln n}{2}-1\rt\},\]
and observe analogously to~\eqref{n(x)capl} that $R_x\subseteq R$ and
\begin{align}
\label{n(x)capr}
\Pr \lt[\lt(|N_G(x)\cap R_x|\geq \frac{np}{2\ln\ln n}\rt)\,|\, \lt(|R_x|= n\exp\lt(-\Omega\lt(\ln n /(\ln\ln n)^2\rt)\rt)\rt)\rt]=o(1/n).
\end{align}

We denote by $V'$ the set of all vertices from $[n]$
whose in- an outdegrees in $\vG$ lie in $\lt(\frac{np-np/\ln\ln n}{2}, \frac{np+np/\ln\ln n}{2}\rt)$.
Notice that $[n]\setminus V'\subseteq L_x\cup R_x$ for every $x\in [n]$.
Hence, we see that $|V'|>n-\frac{n}{\ln^2 n}$ a.a.s. by~\eqref{|L_x|} and~\eqref{|R_x|}.
Furthermore, from~\eqref{n(x)capl} and~\eqref{n(x)capr}
we obtain that all in- an outdegrees in $\vG[V']$ lie in $\lt(\frac{np-3np/\ln\ln n}{2}, \frac{np+np/\ln\ln n}{2}\rt)$ a.a.s.,
completing the proof of the lemma.
\proofend

From now on, whenever we have $n$ and $p$ chosen, we denote
$$
d=d(n,p)=np-100np/\ln\ln n\,.
$$

In the following lemma, we show that the random graph contains
a.a.s. many $2$-factors.

\begin{lemma}
\label{number2-factors}
The random graph $G\sim G(n,p)$ with $p\geq \ln n/n$ satisfies
\[ \sum  _{s\in [n/3]}2^s f'(G,s)\geq d^{-n/\ln^2 n}(d/e)^n\]
a.a.s.
\end{lemma}

\proofstart
In order to use Lemma~\ref{orient},
we orient $G$ at random to obtain $\vG$
(as always in this paper, for every edge each of the two possible orientations gets probability $1/2$ independently of the choices of other edges).

First we show that the second condition of Lemma~\ref{orient} holds a.a.s. for $\vG$ with the intuitive choice $r=np/2$.
Since the maximum degree of $G$ is at most $3np$ a.a.s. (see e.g.~\cite{bollobasbook}),
we obtain that in $G$ a.a.s. for any two sets $A$ and $B$  with $|A|>100|B|$ the number of edges between them is at most
$3np|B|<0.4np\sqrt{|A||B|}$.
Hence, $e_{\vG}(A,B)\leq 0.4np\sqrt{|A||B|}$ and $e_{\vG}(B,A)\leq 0.4np\sqrt{|A||B|}$.
Thus, we are left with the case of sets $A$ and $B$ of sizes $|A|\leq 100|B|$ and $|B|\leq 100 |A|$.

For small disjoint sets, we obtain using Chernoff's bound (see e.g. Theorem~A.1.12 in~\cite{alonspencer})
\begin{align*}
&\Pr\lt(\exists A',B'\subset [n], \, A'\cap B'=\emptyset,\, |A'||B'|\leq \frac{n^2}{\ln\ln n},\, |A|\leq 100|B|\leq 10^4 |A|:\, e_{\vG}(A',B')\geq 0.4np\sqrt{|A'||B'|}\rt)\\
&\leq
\sum_{a,b=o(n),\, a=\Theta(b)}\binom{n}{a}\binom{n-a}{b}\exp\lt(-\Omega\lt(np\sqrt{ab}\ln\lt(\frac{np\sqrt{ab}}{pab}\rt)\rt)\rt)\\
&\leq
\sum_{a,b=o(n)}\lt(\frac{ne}{a}\rt)^{a}\lt(\frac{ne}{b}\rt)^{b}\exp\lt(-\Omega\lt(a\ln n\ln\lt(\Omega\lt(\frac{n}{a}\rt)\rt)\rt)
-\Omega\lt(b\ln n\ln\lt(\Omega\lt(\frac{n}{b}\rt)\rt)\rt)\rt)\\
&\leq \sum_{a,b=o(n)}\exp\lt(a\ln\lt(\frac{ne}{a}\rt)+b\ln\lt(\frac{ne}{b}\rt)
-\Omega\lt(a\ln\lt(\frac{n}{a}\rt)\ln n\rt)-\Omega\lt(b\ln\lt(\frac{n}{b}\rt)\ln n\rt)\rt)\\
&=o(1).
\end{align*}
Similarly, for large disjoint sets we obtain using Chernoff's bound (see e.g. Corollary~A.1.14 in~\cite{alonspencer})
\begin{align*}
&\Pr\lt(\exists A',B'\subset [n], \, A'\cap B'=\emptyset,\, |A'||B'|> \frac{n^2}{\ln\ln n},\, |A|,|B|\leq 0.6n:\, e_{\vG}(A',B')\geq 0.4np\sqrt{|A'||B'|}\rt)\\
&\leq
\sum_{a,b\leq n,\, ab> \frac{n^2}{\ln\ln n} }\binom{n}{a}\binom{n-a}{b}\exp\lt(-\Omega\lt(abp\rt)\rt)\\
&\leq 4^n \exp\lt(-\Omega\lt(\frac{n\ln n}{\ln\ln n}\rt)\rt)
=o(1).
\end{align*}
Hence, a.a.s. for every pair of disjoint sets $A'$ and $B'$,
the number of edges going from $A'$ to $B'$ satisfies
\begin{equation}
\label{e(a,b)disjoint}
 e_{\vec{G}}(A',B')<0.4np\sqrt{|A'||B'|}.
\end{equation}

Analogously, we see that a.a.s. for every $M\subseteq [n]$ of size at most $0.6 n$,
\begin{equation}
\label{e(m)}
e_G(M)<0.4np|M|.
\end{equation}

Thus, a.a.s. for every $A,B\subset [n]$ of size $|A|,|B|\leq 0.6 n$,
the number of edges going from $A$ to $B$ in $\vG$ is bounded by
\begin{align*}
 e_{\vG}(A,B)=e_ {\vG}(A\setminus B,B\setminus A)+e_G(A\cap B)
&\overset{\eqref{e(a,b)disjoint},\, \eqref{e(m)}}{<}
0.4np\sqrt{|A\setminus B||B\setminus A|}+0.4np|A\cap B|\\
&\leq 0.4np\sqrt{|A||B|},
\end{align*}
establishing that the second condition of Lemma~\ref{orient} holds a.a.s. for every subgraph of $\vG$.

Hence, by Lemma~\ref{core} the graph $G\sim G(n,p)$ a.a.s. is such that for a random orientation $\vG$,
there a.a.s. exists a vertex set $V'\subseteq [n]$ of size at least $n-n/\ln^2 n$
such that the induced subgraph $\vG[V']$ satisfies the conditions of Lemma~\ref{orient} with $r=np/2$.
Applying Lemma~\ref{orient} to this induced subgraph, we obtain
\[\sum  _{s\in [n/3]} f\lt(\vG[V'],s\rt)\geq \lt(\frac{d}{2e}\rt)^{n-n/\ln^2n}\] a.a.s.
Thus, we obtain
\[\sum_{s\in [n/3]} \Exp\lt(f'\lt(\vG,s\rt)\rt)\geq (1-o(1))\lt(\frac{d}{2e}\rt)^{n-n/\ln^2n}\] a.a.s.,
where the expectation is taken over the random choice of orienting the edges of $G$, the process creating $\vec{G}$ from $G$.

On the other hand, when we orient the edges, an almost $2$-factor of $G$ with exactly $s$ cycles becomes an almost $1$-factor of $\vG$
with probability at most $2^{\frac{n}{\ln^2n}-n+s}$,
implying
\[\sum  _{s\in [n/3]}2^s f'(G,s)\geq \sum_{s\in [n/3]}2^{n-\frac{n}{\ln^2n}} \Exp\lt(f'\lt(\vG,s\rt)\rt).\]

Putting these two facts together, we obtain
\[\sum  _{s\in [n/3]}2^s f'(G,s)\geq  (1+o(1))\lt(\frac{d}{e}\rt)^{n-n/\ln^2n}\geq d^{-n/\ln^2 n}(d/e)^n\] a.a.s., completing the proof of the lemma.
\proofend

We show now that there are typically many almost $2$-factors in $G$ with a small number of cycles.
We denote \[s^{*}= s^{*}(n)=\frac{n}{\ln n\sqrt{\ln\ln n}}.\]

\begin{lemma}
\label{shortcycles}
For every $p\geq \ln n/n$, the random graph $G\sim G(n,p)$ satisfies
\[\sum_{s= 1}^{s^{*}}f'(G, s)\geq\left(np/e\right)^n (1-o(1))^n\]
a.a.s.
\end{lemma}

\proofstart
By Lemma~\ref{number2-factors} we know that $\sum  _{s\in [n/3]} 2^s f'(G,s)\geq d^{-n/\ln^2 n}(d/e)^n$ a.a.s.

We show now that the contribution of almost 2-factors with too many cycles is negligible.
We use the estimate~(5) of~\cite{KKO1}:
in the random graph $H\sim G(n',p)$, for every $s\geq \ln n'$,
\[\E(f(H,s))\leq \frac{(n'-1)!\lt(\ln n'\rt)^{s-1}p^{n'}}{(s-1)!2^s}.\]
We obtain
\begin{align*}
\sum_{s= s^{*}}^{n/3}\E(2^sf'(G, s))
&\leq\sum_{\ell \leq n/\ln^2 n}\binom{n}{\ell}\sum_{s= s^{*}}^{n/3}\frac{n! \lt(\ln n\rt)^s p^{n-\ell}}{s!}\\
&\leq n!p^n \lt(\frac{n}{p}\rt)^{n/\ln^2 n}\sum_{s= s^{*}}^{n/3}\lt(\frac{s}{e\ln n}\rt)^{-s}\\
&= (d/e)^ne^{O(n/\ln\ln n)} \lt(\frac{s^{*}}{e\ln n}\rt)^{-s^{*}}\\
&=o\left(d^{-n/\ln^2 n}\left(d/e\right)^n\right).
\end{align*}

Hence, using this estimate together with Markov's inequality,
we see that the number of almost $2$-factors of $G$ with at most $s^{*}$ cycles is
\[\sum_{s=1}^{ s^{*}}f'(G, s)\geq \frac{1}{2}2^{-s^{*}}d^{-n/\ln^2 n}\left(d/e\right)^n
=\left(np/e\right)^n(1-o(1))^n\]
a.a.s.
\proofend

The next technicality we need to prove in order to be ready to prove the main theorem is the expansion of $G(n,p)$.

To collect all but one expansion properties that we need, we make the following definition.
\begin{defi}
\label{def}
We call a graph $G$ with the vertex set $[n]$ a $p$-expander, if there exists a set $D\subset [n]$
such that $G$ and $D$ satisfy the following properties:
\begin{itemize}
\item
$|D|\leq n^{0.09}.$
\item
The graph $G$ does not contain a non-empty path of length at most $\frac{2\ln n}{3\ln\ln n}$
such that both of its (possibly identical)
endpoints lie in $D$.
\item
For every set $S\subset [n]\setminus D$ of size $|S|\leq \frac{1}{p}$,
its external neighborhood satisfies $|N(S)|\geq \frac{np}{1000}|S|$.
\end{itemize}
\end{defi}

The following lemma shows that these properties are pseudo-random.
\begin{lemma}
\label{expander}
Consider the two-round expansion of the random graph and fix $G\sim G(n,p)$ with $\ln n/n\leq p\leq 1-2\ln\ln n/n$
and $G\subseteq \hat{G}\sim G(n,\hat{p})$ with $\hat{p}=p+2\ln\ln n/n$.
Then it is a.a.s. true that every graph $G'$ satisfying $G\subseteq G'\subseteq \hat{G}$ is a $p$-expander.
\end{lemma}

\proofstart
We first expose $G$ and fix $D=\{v\in [n]:\, d_G(v)<np/100\}$ to be the set of all vertices of $G$ with degree less than $np/100$ in $G$.
Since for a fixed set $D$ the second property is decreasing
and the third property is increasing,
it suffices to prove the second statement for $\hat{G}$
and the third statement for $G$.

The first property is satisfied by Claim~4.3 from~\cite{BSKS} a.a.s.
The second property can be proven to hold in $\hat{G}$ a.a.s. similarly to Claim~4.4 from~\cite{BSKS}
(there it is proven to hold for $G$ a.a.s.)

For the third property, assume to the contrary that there exists a set $S\subset [n]\setminus D$ of size at most $|S|\leq \frac{1}{p}$
such that its external neighborhood in $G$ satisfies $|N_G(S)|< \frac{np}{1000}|S|$.
By the definition of $D$, the number of edges incident to $S$ in $G$ is
\[ e_G(S, N_G(S)\cup S)\geq |S|np/200.\]
But Chernoff's inequality (see e.g Theorem~A.1.12 in~\cite{alonspencer}) tells us that
\begin{align*}
&\Pr\lt(\exists A,B\subseteq [n], \, |A|\leq \frac{1}{p}, |B|<\frac{np}{1000}|A|:\, e_G(A, B\cup A)\geq |A|np/200\rt)\\
&< \sum_{A,B\subset [n],\, |A|\leq 1/p, \, |B|<|A|np/1000}\binom{n}{|A|}\binom{n-|A|}{|B|}
\lt(e\cdot \frac{\E\lt(\lt|e_G(A, B\cup A)\rt|\rt)}{|A|np/200}\rt)^{|A|np/200}\\
&< \sum_{A,B\subset [n],\, |A|\leq 1/p, \, |B|<|A|np/1000}\binom{n}{|A|}\binom{n}{|B|}\lt(\frac{200e|A||A\cup B|p}{|A|np}\rt)^{|A|np/200}\\
&< \sum_{a\leq 1/p}anp\binom{n}{a}\binom{n}{\frac{anp}{1000}}\lt(\frac{3ap}{5}\rt)^{anp/200}\\
&< \sum_{a\leq 1/p}anp\lt(\frac{3ap}{5}\rt)^{anp/400}
=o(1),\end{align*}
providing that the third property holds in $G$ a.a.s.
\proofend

The proof of the next lemma is based on the ingenious
rotation-extension technique, developed by P\'osa~\cite{Posa}, and
applied later in a multitude of papers on Hamiltonicity, mostly of
random or pseudorandom graphs (see for example \cite{BFF1},~\cite{FK},~\cite{KS},~\cite{KrSu}).

Let $G$ be a graph and let $P_0=(v_1,v_2,\ldots,v_q)$ be a path in $G$.
If $1 \leq i \leq q-2$ and $(v_q,v_i)$ is an edge of $G$, then there exists a path
$P'=(v_1 v_2\ldots v_i v_q v_{q-1} \ldots v_{i+1})$ in $G$ with the same set of vertices.
The path $P'$ is called a {\em rotation} of $P_0$ with {\em
fixed endpoint} $v_1$ and {\em pivot} $v_i$. The edge
$(v_i,v_{i+1})$ is called the {\em broken} edge of the rotation. We
say that the segment $v_{i+1} \ldots v_q$ of $P_0$ is {\em reversed} in
$P'$.
In case the new endpoint $v_{i+1}$ has a neighbor $v_j$ such that
$j \notin \{i, i+2\}$, then we can rotate $P'$ further to obtain
more paths of the same length.
We will use rotations
together with the expansion properties from Lemma~\ref{expander} and the necessary minimum degree condition
to find a path on the same vertex set as $P_0$ with
large rotation endpoint sets.

The next lemma shows that in any graph having the $p$-expander property and minimum degree $2$,
for any path $P_0$ and its endpoint $v_1$,
after a small number of rotations with fixed endpoint $v_1$,
we either create many other endpoints or extend the path.
Its proof has certain similarities to the proofs of Lemma~8 from~\cite{us+T} and of Claim~2.2 from~\cite{ProfDrSzabo}.

\begin{lemma}
\label{St}
Let $n$ be a sufficiently large integer and $G$ be an $n$-vertex $p$-expander with minimum degree $\d(G)\geq 2$ and $np\geq \ln n$.
Let $P_0$ be a $v_1w$-path in $G$.
Denote by $B(v_1)\subset
V(P_0)$ the set of all vertices $v\in
V$ for which there is a $v_1v$-path on the vertex set $V(P_0)$ which can be
obtained from $P_0$ by at most $3\frac{\ln n}{\ln (np)}$ rotations
with fixed endpoint $v_1$.
Then $B(v_1)$ satisfies one of the following properties:
\begin{itemize}
\item
there exists a vertex $v\in B(v_1)$ with a neighbor outside $V(P_0)$, or
\item
$|B(v_1)|\geq n/3000$.
\end{itemize}
\end{lemma}

\proofstart
Assume that $B(v_1)$ does not have the first property
(i.e., for every $v\in B(v_1)$ it holds that $N(v)\subseteq V(P_0)$).

Let $t_0$ be the smallest integer such that
$\left(\frac{np}{3000}\right)^{t_0-1} \geq \frac{1}{p}$; note that $t_0
\leq  2\frac{\ln n}{\ln (np)}$.
Since $G$ is a $p$-expander, there is a corresponding vertex set $D$ as in Definition~\ref{def}.

At the first step,
we find a neighbor $u\not \in D\cup N(D)$ of $w$ that is not a neighbor of $w$ along $P_0$.
Its existence is guaranteed since
$w$ has at least two neighbors along $P_0$,
and by the second $p$-expansion property,
at most one of them can have a neighbor in $D$.
We rotate the initial path $P_0$ with pivot $u$
and call the resulting path $P'=(v_1, \ldots, v_q)$.
Notice that this way, $v_q$ is guaranteed not to belong to $D$.

We construct a sequence of sets $S_0, \ldots , S_{t_0}
\subseteq B(v_1)\setminus D \subseteq V(P_0)\setminus\{v_1\}$ of vertices, such that for
every $0 \leq t \leq t_0$ and every $v\in S_t$, $v$ is the endpoint of a
path which can be obtained from $P'$ by a sequence of $t$ rotations with fixed
endpoint
$v_1$, such that for every $0 \leq i < t$,
the non-$v_1$-endpoint of the path after the $i$th rotation
is contained in $S_i$.
Moreover, $|S_t|=\left(\frac{np}{3000}\right)^{t}$ for every $t\leq
t_0-2$, $|S_{t_0-1}|=\frac{1}{p}$, and $|S_{t_0}|\geq\frac{n}{3000}$.

We construct these sets by induction on $t$. For $t=0$, one can
choose $S_0=\{ v_q\}$ and all requirements are trivially satisfied.

Let now $t$ be an integer with $0< t\leq t_0-1$ and assume that the sets $S_0,
\ldots , S_{t-1}$ with the appropriate properties have already been
constructed. We will now construct $S_t$. Let
\[
T= \lt\{v_i\in N(S_{t-1}) : v_{i-1},v_i,v_{i+1}\not\in
\bigcup_{j=0}^{t-1}S_j\cup D\rt\}
\]
be the set of potential pivots for the $t$th rotation, and notice that $T\subset
V(P_0)$ due to our assumption,
since $T\subseteq N(S_{t-1})$ and $S_{t-1}\subseteq B(v_1)$. Assume now
that $v_i\in T$, $y\in S_{t-1}$ and $(v_i,y)\in E(G)$. Then, by the induction
hypothesis, a $v_1
y$-path $Q$ can be obtained from $P'$ by $t-1$ rotations such that
after the $j$th rotation, the non-$v_1$-endpoint is in $S_j$ for
every $0 \leq j \leq t-1$. Each such rotation breaks an edge which is incident with
the new endpoint, obtained in that rotation. Since $v_{i-1}, v_i, v_{i+1}$ are not
endpoints
after any of these $t-1$ rotations, both edges
$(v_{i-1},v_i)$ and
$(v_i,v_{i+1})$ of the original path $P'$ must be unbroken and thus
must be present in $Q$.

Hence, rotating $Q$ with pivot $v_i$ will make either $v_{i-1}$ or
$v_{i+1}$ an endpoint (which of the two, depends on whether the unbroken
segment $v_{i-1} v_i v_{i+1}$ is reversed or not after the first
$t-1$ rotations). Assume without loss of generality that the endpoint is $v_{i-1}$.
We add $v_{i-1}$
to the set $\hat{S}_{t}$ of new endpoints and say that {\em $v_i$
placed $v_{i-1}$ in $\hat{S}_{t}$}. The only other vertex that can
place $v_{i-1}$ in $\hat{S}_{t}$ is $v_{i-2}$ (if it exists).

Observe now that if $t<0.1\ln n/\ln\ln n$,
the distance between any vertex from $S_{t-1}$ and $v_q$
is at most $2t-2<0.2\ln n/\ln\ln n$ by the way the sets were constructed.
Hence, between any two vertices from $N(S_{t-1})\cup N(N(S_{t-1}))$,
there is a path of length at most $0.5 \ln n/\ln\ln n$.
Thus at most one vertex from $D$ can be in $N(S_{t-1})\cup N(N(S_{t-1}))$.
On the other hand, it $t\geq 0.1\ln n/\ln\ln n$,
then $|D|\leq n^{0.09}=o(|S_{t-1}|)=o(|N(S_{t-1})|)$.
Thus, in both cases $|D\cap (N(S_{t-1})\cup N(N(S_{t-1})))|=o(|N(S_{t-1})|)$.

Combining all this information together, we obtain
\begin{align*}
|\hat{S}_{t}| &\geq \frac{1}{2} |T|\\
&\geq
\frac{1}{2} \left(|N(S_{t-1})|-3(1+|S_1|+\ldots+|S_{t-1}|+|D\cap (N(S_{t-1})\cup N(N(S_{t-1})))|\right))\\
&\geq \left(\frac{np}{3000}\right)^{t}.
\end{align*}
Clearly we can
delete arbitrary elements of $\hat{S}_{t}$ to obtain $S_{t}$ of size
$\left(\frac{np}{3000}\right)^{t}$ if $t\leq t_0-2$ and of size
$\frac{1}{p}$ if $t=t_0-1$. So the proof of the induction
step is complete and we have constructed the sets $S_0, \ldots,
S_{t_0-1}$.

To construct $S_{t_0}$ we use the same technique as
above, only the calculations are slightly different.
\begin{align*}
|\hat{S}_{t_0}| &\ge \frac{1}{2} |T|\\
&\ge
\frac{1}{2}
\left(|N(S_{t_0-1})|-3(1+|S_1|+\ldots+|S_{t_0-2}|+|S_{t_0-1}|+|D\cap (N(S_{t-1})\cup N(N(S_{t-1})))|\right))\\
&\geq n/3000.
\end{align*}

The set $S_{t_0}:=\hat{S}_{t_0}$ is by construction a subset of $B(v_1)$,
and the number of rotations needed to make any of its vertices an endpoint of the current path is at most $t_0+1$,
concluding the proof of the
lemma. \proofend

The proof of the following lemma relies on the final part of the proof of Theorem~1 from~\cite{Michael} and uses Lemma~\ref{St}.
It shows that under certain pseudorandom conditions in a graph $G$ for every almost $2$-factor,
after adding few random edges,
there exists a Hamilton cycle within a small Hamming distance from it a.a.s.

\begin{lemma}
\label{digest}
Let $G$ be a connected $n$-vertex $p$-expander with minimum degree $2$ and $S$ be a set of vertices of $G$ of size $|S|=o(n)$
such that there exist at least $n$ non-edges in $G$ not incident to $S$.
Let $F$ be an almost $2$-factor of $G$ with at most $s^{*}$ cycles.
Choose $n$ non-edges $e_1,\ldots, e_n$ of $G$ i.a.r. under the condition that none of them is incident to $S$
and denote by $G'$ the (random) graph obtained from $G$ by turning them into edges.
Then, if it is a.a.s. true that every graph $G\subseteq \hat{G}\subseteq G'$ is a $p$-expander,
then
$G'$  a.a.s. contains a Hamilton cycle $H$
with Hamming distance at most $17s^{*}\ln n/\ln(np)$ from $F$.
\end{lemma}

\proofstart
Fix an arbitrary component $C \subseteq F$.
Since $G$ is connected,
there exists an edge in $G$ connecting a vertex $v\in V(C)$ and $y\not \in V(C)$ -
unless of course $C$ is already Hamiltonian.
We denote by $C'$ the component of $y$ in $F$.
Opening $C$ up by deleting an edge of $C$ incident to $v$ (no need to do so if $C$ is just one isolated vertex),
we get a path $P$. We append the edge $vy$ to $P$,
go through it to $C'$,
and if $C'$ is a cycle, then we open it up by deleting an edge of $C'$
incident to $y$
to get a longer path $P'$ and repeat the argument.
If at some point there are no edges between the endpoints of the current path $P''$ and other components from $F$,
then we can fix one endpoint $x$ of $P''$ and rotate $P''$ using Lemma~\ref{St} to extend it outside
or to obtain a set $B(x)$
of size at least $|B(x)|\geq n/3000$ of potential other endpoints.
For every vertex $z\in B(x)$,
we can rotate the resulting path fixing $z$ as one endpoint
to obtain a set $A(z)$ of size at least $|A(z)|\geq n/3000$ of potential other endpoints
or to extend the path outside.
If the path still cannot
be extended outside
and we can still not close it to a cycle,
we have a set $E'$ of at least $10^{-8}n^2$ non-edges of $G$ not incident to $S$,
so that turning any of them into an edge would close the path to a cycle.
We add pairs $e_1, e_2,\ldots$ to $E(G)$, until one of them falls inside $E'$.
Notice that for every $i\in [n]$, the pair $e_i$ falls into $E'$ with at least some constant positive probability.
This means that considering events ``$e_i\in E'$'',
every event has probability $\Theta(1)$ regardless of the previous events.
Notice that in a successful round, the number of components gets reduced or a Hamilton cycle is created,
since the edge that appeared in $E'$ closed the path into a cycle or extended the path directly.
To reduce the number of components by one, we do at most $\ln n/\ln (np)$ rotations by Lemma~\ref{St},
therefore increasing the Hamming distance from $F$ by at most $4+12\ln n/\ln (np)$.
Since it is enough to have $s^{*}+n/\ln^2n$ successful events to obtain a Hamilton cycle,
the expected number of needed turns of non-edges into edges is at most $O(1)\cdot(s^{*}+n/\ln^2n)=o(n)$.
Hence the $n$ additional edges suffice to create a Hamilton cycle $H$ from $F$ by Markov's inequality a.a.s.,
replacing at most $8\ln n/\ln(np)$ edges for every component of $F$.
Thus, the Hamming distance between $F$ and $H$ is at most $2\cdot 8\frac{\ln n}{\ln(np)}(s^{*}+n/\ln^2 n)\leq 17s^{*}\ln n/\ln(np)$.
\proofend

We are now ready to prove Theorem~\ref{main}.

\proofstart
Notice that only the lower bound is of interest for us.
We expose $G$ in two rounds.

We choose a function $p_1=p_1(n)$ such that $\frac{\ln n+\ln\ln n+\omega(1)}{n}\leq p_1\leq p-\frac{\omega(1)}{n}$,
$p_1 \leq 1 ? 2 \ln \ln n/n$,
and $p_1=(1-o(1))p$.
In the first round, we expose $G_1\sim G(n,p_1)$.
We determine
$D:=\{v\in [n]:\, d_{G_1}(v)<np_1/100$\}.

In the second round, we expose the binomial random graph $G_2$ by including every edge from $K_n\setminus G_1$ into $E(G_2)$
with probability $p_2:=\frac{p-p_1}{1-p_1}$.
Since $(1-p_1)(p_2)=1-p$, we obtain a graph $G:=G_1\cup G_2 \sim G(n,p)$.
We know by Lemma~\ref{expander} that a.a.s. every graph between $G_1$ and $G$
including them both is a $p_1$-expander.
Furthermore, notice that the expected number of edges in $G_2$ is $\binom{n}{2}p_2\geq \binom{n}{2}(1-p_1)=\omega(n)$,
hence a.a.s. at least $n$ additional random edges appeared in the second round of expansion by Markov's inequality.
Since these edges were chosen i.a.r.,
and $G_1$ is a.a.s. connected with minimum degree at least $2$ (see e.g.~\cite{bollobasbook}),
the conditions of Lemma~\ref{digest} are satisfied for $G_1$ and the first $n$ edges exposed in the second round with $S=\emptyset$.

Now, we put all we know together:
\begin{itemize}
\item
By Lemma~\ref{shortcycles} we obtain
$\sum_{s= 1}^{s^{*}}f'(G_1, s)\geq\left(np_1/e\right)^n (1-o(1))^n$ a.a.s.
\item
For every almost $2$-factor $F$ of $G_1$ with at most $s^{*}$ cycles,
there a.a.s. exists a Hamilton cycle in $G$ with Hamming distance at most $k:=17s^{*}\ln n/\ln(np_1)=\frac{17n}{\ln(np_1)\sqrt{\ln\ln n}}=o(n)$ from $F$
by Lemma~\ref{digest}.
\item
On the other hand, for every Hamilton cycle $H$ in $G$,
to obtain an almost $2$-factor of $G_1$ of distance at most $k$ from $H$,
we can first delete at most $k$ edges of $H$,
thus obtaining a collection of at most $k$ paths.
These paths should then be tailored into an almost $2$-factor,
and the choices here are for each of the at most $2k$ endpoints of the paths
to be connected to one of its $\Delta(G_1)$ neighbors in $G_1$
or to stay isolated.
Thus, there are at most $\binom{n}{k}(\Delta(G_1)+1)^{2k}$
almost $2$-factors of $G$ with Hamming distance at most $k$ from $H$.
\item
Hence, by double counting almost $2$-factors of $G$ with at most $s^*$ cycles, we obtain
\[h(G)\geq \frac{\sum_{s= 1}^{s^{*}}f'(G', s)}{\binom{n}{k}(\Delta(G_1)+1)^{2k}}\geq\left(np_1/e\right)^n \frac{(1-o(1))^n}{2^{o(n)} (4\ln(np_1))^{2k}} =\left(np/e\right)^n (1-o(1))^n\]
\end{itemize}
a.a.s.
\proofend

To strengthen the result of Cooper and Frieze~\cite{CF}, we now prove Theorem~\ref{hittingtime}.

\proofstart
The proof goes along the argument of Theorem~\ref{main}, but now we expose the graph in three rounds.
We first expose $G_1\sim G(n,\ln n/n)$ and fix the set $D$ of vertices of degree at most $\ln n/100$.
Notice that similarly to the argument in the proof of Lemma~\ref{expander},
Claim 4.3 from~\cite{BSKS} implies that $|D|\leq n^{0.09}$.
In the second round of exposure, in addition to $G_1$ we expose those edges that are incident to $D$ one by one,
until the minimum degree becomes two;
the resulting graph is called $G'$.
In the third round of exposure we consider the binomial random graph $G_2$ by including every edge of $K_n\setminus G_1$ not incident to $D$
with probability $p_2:=\frac{\ln\ln n}{2n}$.

Let us denote by $G$ the graph obtained by stopping the random graph process at the moment the minimum degree becomes two.
Notice first that since in the random graph process $\d\lt( G\lt(n, \frac{\ln n+0.5 \ln\ln n}{n}\rt)\rt)=1$ a.a.s.,
we obtain $G'\cup G_2\subseteq G$ a.a.s.,
where by the union of two graphs with vertex sets $[n]$ we denote the graph on the same vertex set
where the union is taken over the edge sets.
Furthermore, observe that since in the random graph process $\d\lt( G\lt(n, \frac{\ln n+2 \ln\ln n}{n}\rt)\rt)\geq 2$ a.a.s.,
we obtain $G'\cup G_2\subseteq G\lt(n, \frac{\ln n+2 \ln\ln n}{n}\rt)$ a.a.s.
(The two statements above can be found e.g. in~\cite{bollobasbook}.)
Finally, $G'$ is connected a.a.s. because of the expansion properties and the fact that the edge set between two linearly large sets
is not empty (see e.g.~\cite{ProfDrSzabo}).

Since $p_2=\omega(1/n)$ and $|D|=o(n)$ a.a.s., we obtain $|E(G_2)|=\omega(n)$ a.a.s.
Furthermore, these edges are random under the only conditions of being non-edges of $G_1$ and being not incident to $D$.
Hence,
the conditions of Lemma~\ref{digest} are satisfied for $G'$ and the first $n$ edges exposed in the third round.

Thus, following the lines of the proof of Theorem~\ref{main},
we obtain the desired estimate.
\proofend

\section{Concluding remarks}

In this paper we have proven that for any value of the edge
probability $p=p(n)$, for which the random graph $G~G(n,p)$ is
a.a.s. Hamiltonian, the number of Hamilton cycles in $G$ is
$n!p^n(1+o(1))^n$ a.a.s., thus being asymptotically equal to the
expected value -- up to smaller order exponential terms. Of course,
it would be very nice to extend Janson's result~\cite{Janson} to
smaller values of $p$ and to understand more accurately the
distribution of the number of Hamilton cycles in relatively sparse
random graphs. However, given that the machinery used in
\cite{Janson} is rather involved, and the result (limiting
distribution) is somewhat surprising, this will not necessarily be
an easy task.

Our bound on the number of Hamilton cycles in $G(n,p)$ can be used to bound the number of perfect matching
similarly to~\cite{Michael}. Let $m(G)$ denote the number of perfect matchings in the graph $G$.
Since every Hamilton cycle is a union of two perfect matchings, we obtain $h(G)\leq \binom{m(G)}{2}$.
Hence, for $G\sim G(n,p)$ the a.a.s. lower bound on $h(G)$ from Theorem~\ref{main} provides the a.a.s. lower bound $m(G)\geq (np/e)^{n/2}(1-o(1))^n$.
Since the upper bound is easily obtained from the expected value by Markov's inequality similarly to the first paragraph of Section~\ref{proofs},
we have $m(G)= (np/e)^{n/2}(1-o(1))^n$.
The corresponding hitting time statement is obtained by a straightforward modification of the proof of Theorem~\ref{hittingtime}:
In the random graph process, the edge that makes the graph connected a.a.s. creates $(\ln n/e)^{n/2}(1-o(1))^n$ perfect matchings.

\end{document}